
\input gtmacros       
\input gtmonout
\volumenumber{1}
\volumeyear{1998}
\volumename{The Epstein birthday schrift}
\pagenumbers{159}{166}
\input epsf
\input newinsert
\received{15 November 1997}
\published{22 October 1998}
\papernumber{8}

\def\proc#1{\vskip-\lastskip\par\medskip\bf%
\noindent#1\ \number\sectionnumber.\number\resultnumber
\stdspace\sl\global\advance\resultnumber by 1\ignorespaces}
\def\endproc{\rm\par\medskip} 
\def\proclaim#1{\vskip-\lastskip\par\medskip\bf%
\noindent#1\stdspace\sl\ignorespaces} 

\insertskipamount 4pt plus2pt     
\inserthardskipamount 2pt          
\def\prf{\vskip-\lastskip\par\medskip\noindent{\bf Proof}%
\stdspace\rm} %

\def\la{\langle}
\def\ra{\rangle}
\def\Z{{\Bbb Z}}
\def\inv{^{-1}}
\def\ep{\varepsilon}
\def\ssp{\stdspace}

%
%
\reflist

\def\,{\thinspace}\def\Co{{\rm,\ }}

\refkey\CG {\bf A Clifford\Co R\,Z Goldstein}, {\it Tesselations of $S^2$ and
equations over torsion-free groups}, Proc. Edinburgh Maths. Soc.
38 (1995) 485--493 

\refkey\FR {\bf Roger Fenn\Co Colin Rourke}, {\it Klyachko's methods and the
solution of equations over torsion-free groups}, l'Enseign. Math.
42 (1996) 49--74

\refkey\HNN {\bf G Higman}, {\bf B\,H Neumann}, {\bf Hanna Neumann},
{\it Embedding theorems for groups}, J. London Maths. Soc. 24 (1949) 247--254

\refkey\Kly {\bf A Klyachko}, {\it Funny property of sphere and equations over
groups}, Comm. in Alg. 21 (7) (1993) 2555--2575 

\refkey\Lev {\bf F Levin}, {\it Solutions of equations over groups}, 
Bull. Amer. Math. Soc. 68 (1962) 603--604 

\refkey\Neu {\bf B\,H Neumann}, {\it Adjunction of elements to groups}, 
J. London Math. Soc. 18 (1943) 4--11

\endreflist
%
%
%

\title{Characterisation of a class of equations\\
with solutions over torsion-free groups}   
\shorttitle{Equations over torsion-free groups}   

\author{Roger Fenn\\Colin Rourke}  

\address{School of Mathematical Sciences, Sussex University\\
Brighton, BN1 9QH, UK\\{\rm and}\\
\\Mathematics Institute, University of Warwick\\Coventry, CV4 7AL, UK}

\asciiaddress{School of Mathematical Sciences, Sussex University\\         
Brighton, BN1 9QH, UK\\Mathematics Institute,
University of Warwick\\Coventry, CV4 7AL, UK}

\email{R.A.Fenn@sussex.ac.uk, cpr@maths.warwick.ac.uk}

\abstract
We study equations over torsion-free groups in terms of their
``$t$--shape'' (the occurences of the variable $t$ in the equation).
A $t$--shape is {\sl good\/} if any equation with that shape has a
solution.  It is an outstanding conjecture [\Lev] that all $t$--shapes
are good.  In [\FR] we proved the conjecture for a large class of
$t$--shapes called {\sl amenable}.  In [\CG] Clifford and Goldstein
characterised a class of good $t$--shapes using a transformation on
$t$--shapes called the {\sl Magnus derivative}.  In this note we
introduce an inverse transformation called {\sl blowing up}.
Amenability can be defined using blowing up; moreover the connection
with differentiation gives a useful characterisation and implies that
the class of amenable $t$--shapes is strictly larger than the class
considered by Clifford and Goldstein.
\endabstract

\asciiabstract{%
We study equations over torsion-free groups in terms of their
`t-shape' (the occurences of the variable t in the equation).  A
t-shape is good if any equation with that shape has a solution.  It is
an outstanding conjecture that all t-shapes are good.  In [Klyachko's
methods and the solution of equations over torsion-free groups,
l'Enseign. Maths. 42 (1996) 49--74] we proved the conjecture for a
large class of t-shapes called amenable.  In [Tesselations of S^2 and
equations over torsion-free groups, Proc. Edinburgh Maths. Soc. 38
(1995) 485--493] Clifford and Goldstein characterised a class of good
t-shapes using a transformation on t-shapes called the Magnus
derivative.  In this note we introduce an inverse transformation
called blowing up.  Amenability can be defined using blowing up;
moreover the connection with differentiation gives a useful
characterisation and implies that the class of amenable t-shapes is
strictly larger than the class considered by Clifford and Goldstein.}




\maketitle

\section{Introduction}

Let $G$ be a group. An expression of the form 
$$
\eqalignno{r=g_1t^{\ep_1}g_2t^{\ep_2}g_3\cdots t^{\ep_k}&=1, &(1)\cr}
$$ 
where $k\ge1$, $g_i\in G$ and $\ep=\pm1$, is called an {\sl equation}
over $G$ in the {\sl variable} $t$ with {\sl coefficients}
$g_1,g_2,\ldots ,g_k$. The equation is said to have a {\sl solution}
if $G$ embeds in a group $H$ containing an element $t$ for which (1)
holds. This is equivalent to saying that the natural map
$$
G\longrightarrow {G*\left\langle t\right\rangle \over \left\langle
r=1\right\rangle} 
$$
is injective.

The equation is said to be {\sl reduced} if it contains no subword
$tt\inv$ or $t\inv t$ (ie each coefficient which separates a pair
$t,t\inv$ is non-trivial).  The equation is said to be {\sl cyclically
reduced} if all cyclic permutations are reduced and, unless explicitly
stated otherwise, all equations are assumed to be cyclically reduced.
 
The {\sl $t$--shape} of the word $r$ is the sequence 
$t^{\ep_1}t^{\ep_2}\cdots t^{\ep_k}$.  

We use the abbreviated notation $t^m$ for the sequence $tt\cdots t$
($m$ times) and $t^{-m}$ for the sequence $t^{-1}t^{-1}\cdots t^{-1}$
($m$ times).  We call the $t$--shape $t^m$ ($m\in \Z$, $m\ne0$) a {\sl
power} shape.  If a $t$--shape is not a power then after cyclic
permutation it can be written in the form
$$
t^{r_1}t^{-r_2}t^{r_3}\cdots
t^{-r_u},\; u>1
$$
where each $r_i$ is positive.

The sum $\ep =r_1-r_2+\ldots-r_u$ is called the {\sl degree }of the
$t$--shape. The sum $w=r_1+r_2+\ldots+r_u $ is called the {\sl width}
of the $t$--shape.  Note that the width is the length of the
corresponding equation.

We call a cyclic $t$--shape {\sl good} if any corresponding equation
with torsion-free coefficients has a solution. 
\medskip
{\bf Conjecture}\qua[\Lev]\ssp{\sl All $t$--shapes are good.}
\medskip

The conjecture is a special case of the adjunction problem [\Neu] and
for a brief history, see the introduction to [\FR].  The torsion-free
condition is necessary because the $t$--shape $tt^{-1}$ is good [\HNN] 
but for example the equation $ata^2t\inv=1$ has 
no solution over a group in which $a$ has order 4.

The conjecture is known to be true in many cases.  Levin [\Lev] has
proved that power shapes are good (without the torsion-free
hypothesis).  Klyachko [\Kly] has proved that $t$--shapes of degree
$\pm1$ are good.  Furthermore both Clifford and Goldstein [\CG] and
ourselves [\FR] have extended Klyachko's results to larger classes of
$t$--shapes.  The class of good $t$--shapes in [\CG] are characterised
in terms of the {\sl Magnus derivative} and for definitiveness we will
call them {\sl CG--good}.  The class of good $t$--shapes in [\FR] are
called {\sl amenable}.  No usable characterisation of amenability was
given in [\FR] and it is the purpose of this note to supply such a
characterisation and to compare the two classes.

The rest of the paper is organised as follows.  In the next section
(section 2) we review the Magnus derivative (an
operation on $t$--shapes which we refer to simply as {\sl
differentiation}) and define the class of CG--good shapes.  In section
3 we define another operation on $t$--shapes called {\sl blowing up}
and prove that it is the inverse of differentiation.  Finally in
section 4 we give two simple characterisations of amenable shapes.
The first in terms of blowing up and the second, similar to the
characterisation of CG--good shapes, in terms of differentiation.  We
conclude that the class of amenable shapes is strictly larger than the
class of CG--good shapes.

{\bf Acknowledgements}\ssp We are grateful to Martin Edjvet for
suggesting that there might be a connection between the results of the
Clifford--Goldstein paper and ours.  We thank the referee for helpful
comments.

\section{The Magnus derivative}

Let $T=t^{\ep_1}t^{\ep_2}\cdots t^{\ep_w}$, where $\ep_i=\pm1$, be a
$t$--shape.  We regard $T$ as a cyclic $t$--shape and we define the
cyclic $t$--shape $D(T)$, the {\sl Magnus derivative} or simply {\sl
derivative} of $T$, as follows.

Arrange the signs of the exponent powers around a circle. The
$t$--shape is well defined by this up to cyclic symmetry.  Between each
occurence of $+,+$ insert a new $+$, between each occurence of $-,-$
insert a new $-$ and in all other cases do nothing.  Now delete the original
signs. The remaining cyclic sequence of signs defines a new
$t$--shape, $D(T)$.

For example\qquad
$tttt^{-1}tt^{-1}t^{-1}t\quad{\buildrel D \over\rightarrow}\quad
ttt^{-1}t\quad{\buildrel D \over\rightarrow}\quad tt$.

The following is easy to prove.

\proclaim{Lemma}Let the cyclic $t$--shape $T$ have degree $\varepsilon(T)$ and 
width $w(T)$ then:
\items
\item{\bf1\rm)}$\varepsilon (DT)=\varepsilon (T).$

\item{\bf2\rm)}$w(DT)\leq w(T)$ with equality if and only if $T$ is empty
or a power shape.

\item{\bf3\rm)}$D(T)=T$ if and only if $T$ is empty or a power shape.

\item{\bf4\rm)}$D^\alpha (T)$ is empty or a power shape if $\alpha >w(T)/2$.

\item{\bf5\rm)}If $T=t^{r_1}t^{-r_2}t^{r_3}
\cdots t^{-r_k}$, where $r_i\ge1$, is not a power shape then\nl 
$DT=t^{r_1-1}t^{-r_2+1}\cdots t^{-r_k+1}$.
\qed
\enditems\endproc

We can illustrate the effect of differentiation by looking at the {\sl
graph} of the $t$--shape $T=t^{\ep_1}t^{\ep_2}\cdots t^{\ep_w}$.

This is a function $f=f_T\co[0,w]\rightarrow \re$ defined as follows. 
Define $f(0)=0$ and for integers $i$ in the range $0< i\leq w$
$f(i)=\ep_1+\ep_2+\ldots+\ep_i$. Extend $f$ over the whole 
interval by piecewise-linear interpolation. Notice that the graph of the
$t$--shape starts at $(0,0)$ and finishes at $(w,\varepsilon )$.

Figure \figkey\Diff\ shows the graph of the example above and the effect of
differentiation which `smooths off' the peaks and troughs until a
straight line graph is left.

\fig{\Diff: Differentiation}
\epsfbox{diff.eps}
\endfig

A {\sl clump} in a cyclic $t$--shape is defined to be a maximal connected
subsequence of the form $t^m$ where $|m|>1$.   A {\sl one-clump shape}
is a shape with just one clump, which is not the whole sequence,
ie, after possible cyclic permutation and inversion, 
a shape of the form $t^mt^{-1}(tt^{-1})^r$ where $m>1$ 
and $r\ge0$.  We can now define CG--good.  A $t$--shape is {\sl 
CG--good} if, after a (possibly empty) sequence of differentiations 
it becomes a one-clump shape.

\proclaim{Theorem}{\rm(Clifford--Goldstein [\CG])}\ssp 
All CG--good shapes are good. \qed

\section{Blowing up}

We shall now introduce the notion of {\sl blowing up} of a $t$--shape
which was implicit in [\FR].  

We consider non-cyclic $t$--shapes whose graphs start and end at
level 0 and which lie between levels $-m$ and 0.  Such a $t$--shape
will be called an {\sl $m$--block}.  An $m$--block whose graph reaches
level $-m$ at some point will be called a {\sl full $m$--block}. 

\rk{Definition}{\sl $m$--blow up}\ssp  Start with a given cyclic
$t$--shape.   Between each pair $t\inv t$ (ie at local minima of the
graph) insert a full $m$--block.  Between other pairs
insert a general $m$--block (see figure \figkey\blowup). \ppar

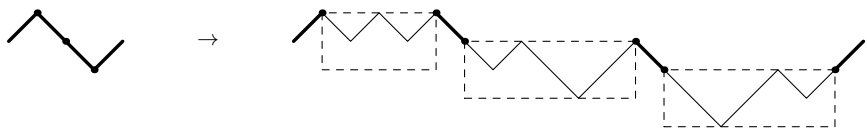
\fig{\blowup: An example of a 2--blow-up}
\epsfxsize 4.5truein
\epsfbox{blow.eps}
\endfig

The definition of blow up is not explicit in [\FR].  
However we shall see later that it coincides with the concept of 
normal form given on page 69 of [\FR].

Notice that a 0--blow up of a shape $T$ is the original shape $T$ but
that, in general, the result of blowing up depends on the choices of the
blocks.  We use the notation $B^m(T)$ for the set of $m$--blow ups of
$T$ and we abbreviate $B^1$ to $B$.

We now prove that blowing up is anti-differentiation.

\proc{Lemma}$U\in B(T)$ if and only if $D(U)=T$.\endproc  

\prf We give a graphical description of $D$.  Start with the graph of
a $t$--shape $T$. Introduce a new vertex halfway along each edge of
the graph.  At each local maximum (respectively minimum) join the new
vertices just below (respectively above) and truncate.  Now contract
the horizontal edges and discard the old vertices.  The result is the
graph of $D(T)$.

This process is illustrated in figure \figkey\Refa, where the new
vertices are open dots and the old vertices are black dots.

\fig{\Refa: Graphical differentiation} 
\epsfxsize 4.7truecm\epsfbox{ref1.eps}
\raise 7pt\hbox{$\rightarrow$}
\epsfxsize 4.7truecm\epsfbox{ref2.eps}
\raise 7pt\hbox{$\rightarrow$}
\epsfxsize 2.35truecm\epsfbox{ref3.eps}
\endfig

To see the connection with 1--blow ups consider the following
alternative description.  Introduce the new vertices
as before but slide them up to the top of the edges.  Discard all the
locally minimal vertices of the graph of $T$ and again reduce the
resulting graph by contracting horizontal edges (see
figure \figkey\Refb).  In this description it is clear that the
discarded pieces are precisely 1--blocks and the lemma follows.
\endprf

\fig{\Refb: Differentiation and 1--blow up} 
\epsfxsize 4.7truecm\epsfbox{ref1.eps}
\raise 7pt\hbox{$\rightarrow$}
\epsfxsize 4.7truecm\epsfbox{ref4.eps}
\raise 7pt\hbox{$\rightarrow$}
\epsfxsize 2.35truecm\epsfbox{ref3.eps}
\endfig

For the next lemma we need to extend differentiation and
blowing up to $m$--blocks.  If $T$ is an $m$--block then we define
an $n$--blow up by inserting full $n$--blocks at local minima and
general $n$--blocks at all other vertices, including the first
and last vertex (in other words we prefix and append a general
$n$--block).  It can then be seen that the $n$--blow up of 
an $m$--block is an $(m+n)$--block and if the original block is
full, then the blow up is also full.
 
We extend differentiation by using the same rule
as for cyclic $t$--shapes.  In graphical terms it has the same
meaning as in the last proof:  Discard all the locally minimal vertices
of the graph and reduce by contracting horizontal edges.  The proof
of the previous lemma then shows that $B$ and $D$ are inverse operations on
$m$--blocks.

\proc{Lemma}{\rm(a)}\ssp $B\circ B^{m}\subset B^{m+1}$\qquad
{\rm (b)}\ssp $DB^{m+1}\subset B^m$.\endproc 

\prf  A 1--blow up of an $m$--blow up can be obtained by 1--blowing 
up the inserted $m$--blocks.  Part (a) now follows from the remarks 
above.  To see part (b) observe that $D$ of a $(m+1)$--blow up
is obtained by differentiating the inserted pieces and thus results
in an $m$--blow up.
\endprf

\proc{Corollary}{\rm(a)}\ssp $B\circ B^{m}=B^{m+1}$\qquad
{\rm(b)}\ssp $B^n=B\circ\ldots\circ B$ ($n$ factors)\nl
{\rm(c)}\ssp $B^n\circ B^m=B^{n+m}$.\endproc 

\prf (a)\ssp
By part (a) of lemma 3.2 we just have to show that if $U\in B^{m+1}(T)$ 
then $U\in B\circ B^{m}(T)$.  But $D(U)\in B^m(T)$ by part (b),  
and $U\in B(D(U))$ by lemma 3.1 and hence $U\in B(D(U))\subset 
B\circ B^{m}(T)$.

Parts (b) and (c) follow by induction.\endprf

\proc{Corollary}$U\in B^n(T)$ if and only if $D^n(U)=T$.\endproc 

\prf Repeat lemma 3.1 $n$ times. \endprf

We now turn to the connection of blowing up with the concept of
normal form defined in [\FR].  

On page 69 of [\FR] we define a word in {\sl normal form} based on a
particular cyclic $t$--shape $T$ as a word obtained from $T$ by
inserting elements of certain subsets ($X$, $J$ and $Y$ defined on
page 65) of the kernel of the exponential map $\ep\co G*\la t\ra
\to\Z$ at {\sl top} (between $t$ and $t\inv$), {\sl middle} (between
$t$ and $t$ or $t\inv$ and $t\inv$) and {\sl bottom} (between $t\inv$
and $t$) positions respectively.  Inspecting the definitions of $X$,
$J$ and $Y$, it can be seen that this corresponds to inserting
$m$--blocks and then allowing a controlled amount of cancellation.  To
be precise, define a {\sl leading string} of an $m$--block to be an
initial string $t\inv t\inv\ldots t\inv$ and a {\sl trailing string}
to be a final string $tt\ldots t$.  Cancellation is allowed for
specified leading and trailing strings of all blocks.  The defining
condition on $X$ is that the graph of the corresponding block must
meet level $0$ after deletion of leading and trailing strings and the
defining condition for $Y$ is that the block must be full.  There is
no condition on $J$.  We call the blocks corresponding to elements of
$X$, $J$ and $Y$, {\sl top}, {\sl middle} and {\sl bottom} blocks,
respectively and we denote the set of words in normal form based on
the cyclic $t$--shape $T$ by $NF(T)$.

\proc{Lemma}$NF(T)=B^m(T)$.\endproc

\prf  Blowing up corresponds to normal form with no cancellation
allowed and hence $NF(T)\supset B^m(T)$.  For the converse suppose
that $U$ is in normal form based on $T$ and that for a particular top
block $D$ the leading $t\inv$ is allowed to cancel.  Define the
$(m-1)$--block $B$ by $D=t\inv BtC$ (see figure \figkey\Cancel).  Then
figure \Cancel\ makes clear that $U$ can also be obtained by appending
$B$ to the block inserted in the previous place and replacing $D$ by
$C$.  After these substitutions there are fewer allowed cancellations.

\fig{\Cancel: The simplification move}
\epsfxsize\hsize
\epsfbox{cancel.eps}
\endfig

Similar arguments simplify the situation if cancellation takes place
at the end of a top block or at either end of a middle block.  (Notice
that no cancellation can take place at bottom blocks.)  Thus by
repeating simplifications of this type a finite number of times, we
see that $U$ is an $m$--blow up of $T$. \endprf

\section{Amenability}

We now recall the definition of {\sl amenable} $t$--shapes from [\FR].

Recall that a {\sl clump} in a cyclic $t$--shape is a maximal
connected subsequence of the form $t^m$ or $t^{-m}$ where $m>1$.
These are said to have {\sl order} $m$ and $-m$ respectively.  We call
a clump of positive order an {\sl up} clump and a clump of negative
order a {\sl down} clump.  A $t$--shape is said to be {\sl suitable}
if it has exactly one up clump which is not the whole sequence and
possibly some down clumps, or if it has exactly one down clump which
is not the whole sequence and possibly some up clumps.  It follows
that, after a possible cyclic rotation or inversion, a suitable
$t$--shape has the form
$$
t^st^{-r_0}tt^{-r_1}t\ldots tt^{-r_k}
$$
where $s>1$, $k\ge0$ and $r_i\ge1$ for $i=0,\ldots,k$. 

We now define
amenable $t$--shapes.  Using lemma 3.5 above we can rephrase the
definition on page 69 of [\FR] as follows.

\rk{Definition}{\sl Amenable $t$--shapes}\ssp A $t$--shape which is
the $m$--blow up of a suitable $t$--shape is called {\sl amenable}.

\proclaim{Theorem}{\rm(Fenn--Rourke [\FR])}\ssp Amenable shapes are good.
\qed\endproc

We now turn to the characterisation of amenability.  Using corollary
3.4, the definition of amenability says that a shape is amenable if
and only if it eventually differentiates to a suitable shape.  But now
a suitable $t$--shape is either a one clump shape or differentiates to
$t^st^{-r}$ for some $r,s\ge 1$.  This in turn either eventually
differentiates to $tt\inv$ or to $t^st\inv$ or to $tt^{-r}$ for some
$r,s\ge 2$.  Now the last two are one clump shapes and so we can see
that a suitable shape either eventually differentiates to a one clump
shape or to $tt\inv$.  To make the final characterisation of
amenability as simple as possible, we make the shape $tt\inv$ an
honorary amenable shape (it is good [\HNN]) and then we have the
following simple characterisation.

\proc{Theorem}{\rm(Characterisation of amenability)}\ssp  A shape is
amenable if and only if, after a (possibly empty) sequence of 
differentiations, it becomes either a one-clump shape or the
shape $tt\inv$. \qed\endproc

\proc{Corollary}  Amenable shapes are a strictly larger class than
CG--good shapes. \qed\endproc

{\bf Final remarks}\qua(1)\ssp The class of amenable shapes which are not
CG--good are precisely those which eventually differentiate to
$tt\inv$: an example would be $tt\inv t^2t^{-2}$.  It seems that the
methods of Clifford and Goldstein can be extended with little extra
work to the smaller class of shapes which eventually differentiate to
the shape $t^2t^{-2}$.  However we cannot see how to extend their
methods to cover all amenable shapes.

(2)\ssp The remark at the top of page 70 of [\FR], which was left
unproven, can be quickly proved using theorem 4.1.

\references         

\Addresses\recd
\bye